\documentclass[12pt]{amsart}

\usepackage{amscd}
\usepackage{amsthm}
\usepackage{amsmath}
\usepackage{latexsym}
\usepackage{graphicx}
\usepackage{delarray}
\usepackage{float}
\usepackage{mathrsfs}
\usepackage{mathrsfs}
\usepackage[psamsfonts]{amssymb}
\usepackage{vmargin}
\usepackage{amsmath}
\usepackage{epsfig,verbatim}
\usepackage{subfigure}
\usepackage{amsmath,amsfonts,amsthm,amsopn,cite,mathrsfs}
\usepackage{amsfonts,amsthm,amsopn}
\usepackage{float}

\begin{document}


\newtheorem{theorem}{Theorem}[section]
\newtheorem{definition}[theorem]{Definition}
\newtheorem{lemma}[theorem]{Lemma}
\newtheorem{proposition}[theorem]{Proposition}
\newtheorem{corollary}[theorem]{Corollary}
\newtheorem{example}[theorem]{Example}
\newtheorem{remark}[theorem]{Remark}

\def\noo{\noindent }

\newcommand{\abs}[1]{\lvert#1\rvert}
\newcommand{\tr}[1]{\operatorname{tr}(#1)}
\newcommand{\absbig}[1]{\big\lvert#1\big\rvert}
\newcommand{\norm}[1]{\lVert#1\rVert}
\newcommand{\norminf}[1]{\norm{#1}_\infty}
\newcommand{\scp}[1]{\langle#1\rangle}
\newcommand{\spn}[1]{\operatorname{span}(#1)}
\newcommand{\scpbig}[1]{\big\langle#1\big\rangle}
\newcommand{\scpBig}[1]{\Big\langle#1\Big\rangle}
\newcommand{\set}[2]{\big\{ \, #1 \, \big| \, #2 \, \big\}}
\newcommand{\inner}[2]{\langle #1,#2\rangle}
\newcommand{\spec}[1]{\operatorname{Sp}(#1)}
\newcommand{\UU}{\mathbb U}
\newcommand{\RR}{\mathbb{R}}
\newcommand{\RRt}{\mathbb{R}^2}
\newcommand{\RRtt}{\RR\times\widehat\RR}
\newcommand{\R}{\mathbb{R}}
\newcommand{\ZZ}{\mathbb{Z}}
\newcommand{\CC}{\mathbb{C}}
\newcommand{\TT}{\mathbb{T}}
\newcommand{\NN}{\mathbb{N}}
\newcommand{\phihat}{\widehat{\phi}}
\newcommand{\lahat}{\widehat{\la}}
\newcommand{\fhat}{\widehat{f}}
\newcommand{\fkhat}{\widehat{f_k}}
\newcommand{\La}{\Lambda}
\newcommand{\Om}{\Omega}
\newcommand{\Lap}{\La^\perp}
\newcommand{\la}{\lambda}
\newcommand{\dLa}{{\La^\circ}}
\newcommand{\dla}{{\la^\circ}}
\newcommand{\dmu}{{\mu^\circ}}
\newcommand{\al}{\alpha}
\newcommand{\be}{\beta}
\newcommand{\tal}{\tilde{\alpha}}
\newcommand{\tbe}{\tilde{\beta}}
\newcommand{\ga}{\gamma}
\newcommand{\de}{\delta}
\newcommand{\si}{\sigma}
\newcommand{\eps}{\varepsilon}
\newcommand{\ka}{\kappa}
\newcommand{\ZZt}{\ZZ^2}
\newcommand{\ZZd}{\ZZ^d}
\newcommand{\ZZtd}{\ZZ^{2d}}
\newcommand{\klinZ}{k,l\in\ZZ}
\newcommand{\klinZd}{k,l\in\ZZd}
\newcommand{\kiZd}{{k\in\ZZd}}
\newcommand{\laiLa}{\lambda\in\Lambda}
\newcommand{\sumkZn}{\sum_{\kiZn}}
\newcommand{\sumla}{\sum_{\laiLa}}
\newcommand{\sumlap}{\sum_{\laiLap}}
\newcommand{\RRd}{{\RR^d}}
\newcommand{\RRtd}{{\RR^{2d}}}
\newcommand{\FF}{\mathcal{F}}
\newcommand{\M}{\mathcal{M}}
\newcommand{\LL}{\mathcal{L}}
\newcommand{\cH}{\mathcal{H}}
\newcommand{\G}{\mathcal{G}}
\newcommand{\Null}{\mathcal{N}}
\newcommand{\til}{\widetilde}
\newcommand{\lOg}{\ell^1(G)}
\newcommand{\lola}{\ell^1(\Lambda)}
\newcommand{\lolad}{\ell^1(\Lambda^\perp)}
\newcommand{\lolaa}{\ell^1(\Lambda^\circ)}
\newcommand{\lolaac}{\ell^1(\Lambda^\circ ,\ovl c)}
\newcommand{\loab}{\ell^1(\alpha\ZZ\times\beta\ZZ)}
\newcommand{\loba}{\ell^1(\frac{1}{\beta}\ZZ\times\frac{1}{\alpha}\ZZ)}
\newcommand{\LO}{L^1}
\newcommand{\LOR}{\L^1(\RR)}
\newcommand{\LORd}{\LO(\RRd)}
\newcommand{\LOg}{\LO(G)}
\newcommand{\Lt}{L^2}
\newcommand{\LtRd}{\Lt(\RRd)}
\newcommand{\LtRtd}{\Lt(\RRtd)}
\newcommand{\LtR}{\Lt(\RR)}
\newcommand{\LtG}{\Lt(G)}
\newcommand{\id}{\operatorname{id}}
\newcommand{\sseq}{\subseteq}
\newcommand{\supp}{\operatorname{supp}}
\newcommand{\Fro}{\operatorname{Fro}}
\newcommand{\om}{\omega}
\newcommand{\Rarr}{\Rightarrow}
\newcommand{\Lrarr}{\Leftrightarrow}
\newcommand{\ovl}{\overline}
\newcommand{\nob}{\nobreakdash}
\newcommand{\So}{S_0}
\newcommand{\SoR}{\So(\RR)}
\newcommand{\SoRd}{\So(\RRd)}
\newcommand{\SoRt}{\So(\RR^{2d})}
\newcommand{\SoG}{\So(G)}
\newcommand{\Sod}{\So'}
\newcommand{\SodR}{\Sod(\RR)}
\newcommand{\SodRd}{\Sod(\RRd)}

\newcommand{\C}{\mathbf C}
\newcommand{\Z}{\mathbf Z}
\newcommand{\ol}{\overline}
\renewcommand{\phi}{\varphi}
\newcommand{\A}{\mathcal A}
\newcommand{\B}{\mathcal B}
\newcommand{\Bh}{\B(\cH)}
\newcommand{\cC}{\mathcal C}
\newcommand{\K}{\mathcal K}
\newcommand{\LA}{\mathcal L}
\newcommand{\MM}{\mathcal M}
\newcommand{\cG}{\mathcal G}

\newcommand{\CS}{C^*}
\newcommand{\RRdd}{\RRd\times\widehat\RRd}
\newcommand{\GG}{G\times\widehat G}
\newcommand{\pHi}{\varphi}
\newcommand{\piab}{\pi(\alpha k,\beta l)}
\newcommand{\piba}{\pi\big(\frac{k}{\beta},\frac{l}{\alpha}\big)}
\newcommand{\ZZab}{\al\ZZ\times\be\ZZ}
\newcommand{\ZZdab}{\al\ZZd\times\be\ZZd}
\newcommand{\ZZba}{\frac{1}{\be}\ZZ\times\frac{1}{\alpha}\ZZ}
\newcommand{\ZZdba}{\frac{1}{\be}\ZZ^d\times\frac{1}{\alpha}\ZZ^d}
\newcommand{\LOv}{L^1_v}
\newcommand{\LOvs}{L^1_{v_s}}
\newcommand{\LOvRd}{L^1_v(\RRd)} 
\newcommand{\LPqm}{L^{p,q}_m}
\newcommand{\Lqpm}{L^{p',q'}_{1/m}}
\newcommand{\LqpmRt}{L^{p',q'}_{1/m}(\RRtd)}
\newcommand{\Lpqm}{L^{p,q}_m(\RRd)}
\newcommand{\LpqmRt}{L^{p,q}_m(\RRtd)}
\newcommand{\lpqm}{\ell^{p,q}_m}
\newcommand{\lpqtm}{\ell^{p,q}_{\til m}}
\newcommand{\lpqmZtd}{\ell^{p,q}_m(\ZZtd)}
\newcommand{\lpqtmZtd}{\ell^{p,q}_{\til m}(\ZZtd)}

\newcommand{\lOs}{\ell^{1}_{v_s}}
\newcommand{\lOsZ}{\ell^{1}_{v_s}(\ZZ)}
\newcommand{\lOsZd}{\ell^{1}_{v_s}(\ZZd)}
\newcommand{\lOsZtd}{\ell^{1}_{v_s}(\ZZtd)}

\newcommand{\lO}{\ell^{1}}
\newcommand{\lOZ}{\ell^{1}(\ZZ)}
\newcommand{\lOZd}{\ell^{1}(\ZZd)}
\newcommand{\lOZtd}{\ell^{1}(\ZZtd)}
\newcommand{\lOLa}{\ell^{1}(\La)}
\newcommand{\lOdLa}{\ell^{1}(\dLa)}
\newcommand{\Ss}{\mathscr S}
\newcommand{\SsR}{\Ss(\RR)}
\newcommand{\SsRd}{\Ss(\RRd)}
\newcommand{\SsRtd}{\Ss(\RR^{2d})}
\newcommand{\SsRdd}{\Ss(\RRd\times\widehat\RRd)}

\newcommand{\SsLa}{{\mathscr S}(\La)}
\newcommand{\Ssd}{\mathscr {S}'}
\newcommand{\SsdR}{\mathscr {S}'(\RR)}
\newcommand{\SsdRd}{\mathscr {S}'(\RRd)}
\newcommand{\SsdRtd}{\mathscr {S}'(\RR^{2d})}
\newcommand{\SsdRdd}{\mathscr {S}'(\RRd\times\widehat\RRd)}
\newcommand{\MO}{M^{1}}
\newcommand{\MOR}{M^{1}(\RR)}
\newcommand{\MORt}{M^{1}(\RRt)}
\newcommand{\MORd}{M^{1}(\RRd)}
\newcommand{\MORtd}{M^{1}(\RRtd)}
\newcommand{\MOvs}{M^{1}_{v_s}}
\newcommand{\MOvsR}{M^{1}_{v_s}(\RR)}
\newcommand{\MOvsRt}{M^{1}_{v_s}(\RRt)}
\newcommand{\MOvsRd}{M^{1}_{v_s}(\RRd)}
\newcommand{\MOvsRtd}{M^{1}_{v_s}(\RRtd)}
\newcommand{\Mpqm}{M^{p,q}_{m}}
\newcommand{\MpqmR}{M^{p,q}_{m}(\RR)}
\newcommand{\MpqmRt}{M^{p,q}_{m}(\RRt)}
\newcommand{\MpqmRd}{M^{1}_{m}(\RRd)}
\newcommand{\MpqmRtd}{M^{p,q}_{m}(\RRtd)}
\newcommand{\Minfsv}{M^{\infty}_{1/v_s}}
\newcommand{\MinfsvR}{M^{\infty}_{1/v_s}(\RR)}
\newcommand{\MinfsvRt}{M^{\infty}_{1/v_s}(\RRt)}
\newcommand{\MinfsvRd}{M^{\infty}_{1/v_s}(\RRd)}
\newcommand{\MinfsvRtd}{M^{\inf}_{1/v_s}(\RRtd)}

\title{On Spectral Invariance of Non-Commutative Tori}

\author{Franz Luef}
\curraddr{Fakult\"at f\"ur Mathematik\\
Nordbergstrasse 15\\ 1090 Wien\\ Austria}
\email{Franz.Luef@univie.ac.at}
\thanks{The author was supported under the EU-project MEXT-CT-2004-517154.}


\date{January 21, 2006 }



\begin{abstract}
Around 1980 Connes extended the notions of geometry to the
non-commutative setting. Since then {\it non-commutative geometry}
has turned into a very active area of mathematical research. As a
first non-trivial example of a non-commutative manifold Connes
discussed subalgebras of rotation algebras, the so-called {\it
non-commutative tori}. In the last two decades researchers have
unrevealed the relevance of non-commutative tori in a variety of
mathematical and physical fields. In a recent paper we have pointed
out that non-commutative tori appear very naturally in Gabor
analysis. In the present paper we show that Janssen's result on good
window classes in Gabor analysis has already been proved in a
completely different context and in a very disguised form by Connes
in 1980. Our treatment relies on non-commutative analogs of Wiener's
lemma for certain subalgebras of rotation algebras by Gr\"ochenig
and Leinert.
\end{abstract}

\maketitle

\section{Introduction}

In \cite{Gab46} D. Gabor proposed the following method for the
transmission of a speech signal $f$. In modern language he
discretized $f$ into a sequence of bits, i.e. strings of the form
$"0100111"$. A natural way to transmit such strings is to send a
pulse $\phi$ of length $1$ at consequent time intervals and of size
according to the amplitude of $f$, i.e. $f=\sum_{k=1}^m
a_k\phi(t-k)$. Now, a speech signal is a {\it band-limited}
function, i.e. $\supp(\fhat)\subseteq[0,\theta]$ for some finite
real $\theta$, we observe
\begin{equation*}
  \fhat(\om)=\sum_{k=1}^ma_k e^{-2\pi ik\om}\widehat\phi(\om)=\Big(\sum_{k=1}^ma_ke^{-2\pi ik\om}\Big)\widehat\phi(\om).
\end{equation*}
Therefore, the essential support of $f$ and of $\phi$ are equal
which suggests a careful choice of the pulse $\phi$. In practice we
have to transmit more than one signal $f$, e.g. a conversation
between a group of people. Gabor's brilliant idea was to shift each
signal on a different frequency band. More precisely, if
$f_1,...,f_n$ are the band-limited signals we want to transmit, then
he suggested to send $f_l$ on the $l$-th frequency band
\begin{equation*}
  f_l(t)=\sum_{k=1}^ma_{kl}e^{-2\pi il\theta t}\phi(t-k).
\end{equation*}
Therefore, the transmission of all signals $f_1,...,f_n$ corresponds
to
\begin{equation*}
   \sum_{l=1}^nf_l(t)=\sum_{l=1}^n\sum_{k=1}^ma_{kl}e^{-2\pi il\theta t}\phi(t-k).
\end{equation*}
Heisenberg's uncertainty principle implies that $\theta$ has to be
at least greater than $1$ since each signal occupies an area of
size greater than or equal to $1$.
\par
The preceding observations indicates to choose for $\phi$ the Gaussian
$\pHi(t)=e^{-\pi t^2}$, since the Gaussian is well-localized in the time-frequency
plane. Therefore Gabor suggested to decompose an arbitrary signal $f\in\LtRd$ into a
series of time-frequency shifts of a Gaussian over $\ZZtd$:
\begin{equation}\label{e:Gabseries}
  f=\sum_{k,l\in\ZZd}a_{kl}e^{2\pi i l} e^{-\pi (t-k)^2}.
\end{equation}
The last equation is the first example of an {\it atomic
decomposition} which in the last thirty years has led to a new field
of mathematical research with an increasing
literature, see the papers of
Feichtinger and Gr\"ochenig, \cite{FG89a,FG89b}, for the most
general treatment of atomic decompositions. In honour of Gabor's
lasting contribution we call decompositions of type
\eqref{e:Gabseries} {\it Gabor series}. Gabor only gave heuristic
arguments on the convergence of \eqref{e:Gabseries} for $f\in\LtRd$
but a rigorous analysis of \eqref{e:Gabseries} had to wait for
Janssen's contribution \cite{Jan81} of the year 1981. The main
result of Janssen says that the convergence of \eqref{e:Gabseries}
holds only in a weak sense of distributions, while Genossar and
Porat \cite{GP92} showed that the iterative approach may fail for
certain $L^2$-functions. 
 Out of Janssen's paper \cite{Jan81}
emerged the new mathematical field of {\it Gabor analysis}.
\par
Gabor series \eqref{e:Gabseries} are built from  so-called {\it
time-frequency shifts}. Recall that the actions of these operators
on  $f\in L^2(\R^d)$ are given as follows: 
\begin{enumerate}
  \item the {\it translation} operator by
\begin{equation*}
  T_xf(t)=f(t-x),\hspace{2.7cm}x\in\R^d,
\end{equation*}
  \item  the {\it modulation} operator by
\begin{equation*}
  M_{\omega}f(t)=e^{2\pi it\cdot\omega}f(t),\hspace{2.4cm}\omega\in\R^d,
\end{equation*}
  \item {\it time-frequency shifts} by
\begin{equation*}
  \pi(x,\omega)f(t)=M_{\omega}T_xf(t)=e^{2\pi i\omega t}f(t-x),\hspace{0.5cm}(x,\omega)\in\R^{2d}.
\end{equation*}
\end{enumerate}

The time-frequency shifts $(x,\omega,\tau)\mapsto\tau M_{\omega}T_x$
for $(x,\omega)\in\R^{2d}$ and $\tau\in\CC$ with $|\tau|=1$ define
the Schr\"odinger representation of the Heisenberg group,
consequently the time-frequency shifts $\pi(x,\omega)$ for
$(x,\omega)\in\R^{2d}$ are a projective representation of the
time-frequency plane $\R^d\times{\widehat \R}^d$. It is an important
fact that time-frequency shifts satisfy the following composition
law:
\begin{equation}\label{TFcomp}
  \pi(x,\omega)\pi(y,\eta)=e^{-2\pi
  ix\cdot\eta}\pi(x+y,\omega+\eta),
\end{equation}
for $(x,\omega),(y,\eta)$ in the time-frequency plane
$\R^d\times{\widehat \R}^d$.
\par
Therefore, Gabor series \eqref{e:Gabseries} are often written in the
following way:
\begin{equation}\label{e:Gabseries1}
  f=\sum_{k,l\in\ZZd}a_{kl}M_l T_k e^{-\pi t^2}.
\end{equation}
The computation of the coefficient ${\bf a}=(a_{kl})$ of
\eqref{e:Gabseries1} for a given signal $f\in\LtRd$ is one of the
main problems in Gabor analysis. The solution of this problem is a
non-trivial task since the building blocks in $\eqref{e:Gabseries1}$
are non-orthogonal, i.e. $\langle M_{l'} T_{k'} e^{-\pi t^2},M_l T_k
e^{-\pi t^2}\rangle\ne 0$ for $(k,l)\ne(k',l')$. Moreover, the
coefficients ${\bf a}$ are in general not unique.
\par
In section \ref{s:Gabor} we recall the solution of the above problem
in terms of frames due to Daubechies, Grossmann and Meyer,
\cite{DGM86}. Furthermore we give the Janssen representation of a
Gabor frame operator and point out the relevance of modulation
spaces in Gabor analyis. In Section 3 we give a quick review of the
central notions of non-commutative geometry. As application of the
general principles we investigate rotation algebras and
non-commutative tori. Finally we present the link between Gabor
analysis and non-commutative tori. In Section 4 we define spectral
invariant Banach and Fr$\acute{e}$chet algebras. We recall the
fundamental results of Gr\"ochenig and Leinert on window classes of
Gabor frames. We close our discussion with a new approach to the
spectral invariance of non-commutative tori in rotation algebras and
indicate shortly Connes original argument of this important result.

\section{Basics of Gabor analysis}\label{s:Gabor}

Let $g\in L^2(\R^d)$ be a Gabor atom then ${\mathcal
G}(g,\alpha,\beta)=\{\piab g:k,l\in\ZZd\}$ is called a {\it Gabor system}.
Since the Balian-Low principle tells us that it is not possible to
construct (this is in contrast to the situation with wavelets) an
orthonormal basis for $L^2(\R^d)$ of this form, starting e.g.\ from
a Schwartz function $g$, interest in Gabor frames arose. A milestone
was the paper by Daubechies, Grossmann  and Meyer, \cite{DGM86},
where the ``painless use'' of (tight) Gabor frames was suggested.
\par In our case the Gabor frame operator has the following form:
\begin{equation*}
  S_{g,\al,\beta}f=\sum_{\klinZd}\langle f,\piab g\rangle\piab g,\quad  f  \in L^2(\R^d).
\end{equation*}
A Gabor system ${\mathcal G}(g,\al,\beta)$ is called a {\it Gabor
frame} if the Gabor frame operator $S_{g,\al,\beta}$ is invertible,
i.e., if there exist some finite, positive real numbers $A,B$ such
that
\begin{equation*}
  A\cdot{\text I}\le S_{g,\al,\beta}\le B\cdot{\text I}
\end{equation*}
or equivalently,
\begin{equation*}
  A\|f\|^2\le\sum_{\klinZd}\big|\langle f,\piab g\rangle\big|^2\le B\|f\|^2,
\end{equation*}
for all $f$ in $L^2(\R^d)$. Gabor frames ${\mathcal G}(g,\al,\beta)$
allow the following reconstruction formulas for functions $f$ in $\LtRd$
\begin{eqnarray}\label{reconstruct}
    f&=&(S_{g,\al,\beta})^{-1}S_{g,\al,\beta}f=\sum_{\klinZd}\langle
  f,\piab g\rangle\piab(S_{g,\al,\beta})^{-1}g\\
  &=&S_{g,\al,\beta}(S_{g,\al,\beta})^{-1}f=\sum_{\klinZd}\langle
  f,\piab(S_{g,\al,\beta})^{-1}g\rangle\piab g.
\end{eqnarray}
Due to its appearance in the reconstruction formulas
$\gamma_0:=(S_{g,\al,\beta})^{-1}g$ is called the (canonical) {\it
dual Gabor atom}. Another important observation is that the coefficients in the
reconstruction formula \eqref{reconstruct} are {\it not unique} and
therefore there are other {\it dual atoms} $\gamma\in L^2(\R^d)$
with
$$S_{g,\gamma,\al,\beta}f:=\sum_{\klinZd}\langle
f,\piab\gamma\rangle\piab g.$$ Some authors call $(g,\gamma)$ a {\it
dual pair} of Gabor atoms if $S_{g,\gamma,\al,\beta}={\text I}$. In
\cite{WR90} the engineers Raz and Wexler characterized all dual
pairs for finite Gabor system over cyclic groups. In the case of a
Gabor frame $\cG(g,\al,\beta)$ the following is equivalent (up to a
technical condition):
\begin{enumerate}
\item $(g,\gamma)$ is a dual pair on $\LtRd$.
\item $(\al \beta)^{-d}\langle M_{\frac{l}{\al}}T_{\frac{k}{\beta}}\gamma,M_{\frac{l'}{\al}}T_{\frac{k'}{\beta}}g\rangle=\delta_{kk'}\delta_{ll'}$.
\end{enumerate}
Nowadays the preceding equivalent conditions are called the {\it
Wexler-Raz biorthogonality relations}. This important result was obtained
independently and with completely different methods by Ron-Shen,
Janssen and Daubechies-H.J. Landau-Z. Landau in
\cite{RS97,Jan95,DLL95}. In the following we focus on Janssens
approach since it provides an unexpected link to the work of Connes
and Rieffel on non-commutative tori, see \cite{Lu,Lu05}.
\par
Tolmieri and Orr proposed a new method for the calculation of the
frame bounds of a Gabor system $\cG(g,\al,\be)$ in \cite{TO95}.
There method relies essentially on the following identity:
\begin{equation*}
   \sum_{\klinZd} \langle f_1,\piab g_1\rangle{\langle \piab g_2,f_2\rangle}=\frac{1}{(\al\be)^{d}}\sum_{\klinZd}\langle g_2,\piba g_1\rangle{\langle \piba f_1,f_2\rangle},
\end{equation*}
for $f_1,f_2,g_1,g_2$ in Schwartz space $\SsRd$. Due to its great
importance in his approach to the Wexler-Raz biorthogonality
relations Janssen called the previous identity the {\it Fundamental
identity of Gabor analysis} (FIGA). In \cite{DLL95} it goes by the
name {\it Weyl-Heisenberg identity} because many researchers call
$\cG(g,\al,\be)$ a {\it Weyl-Heisenberg frame}. The work of Wexler-Raz related 
the original Gabor system $\cG(g,\al,\be)$ with another Gabor system $\cG(g,\frac{1}{\be},\frac{1}{\al})$ over the
lattice $\ZZdba$. Neither Daubechies and her collaborators nor
Janssen offered an explanation for the specific relation between the
lattices $\ZZdab$ and $\ZZdba$ in their work. In \cite{FK98} Feichtinger and Kozek
were able to reveal a group theoretical interpretation for the
pairing of the two lattices. They observed that the commutant of all
time-frequency shifts $\{\piab:\klinZd\}$ is given by all
time-frequency shifts $\{\piba:\klinZd\}$. Hence they were in the
position to generalize the results in \cite{DLL95,Jan95,RS97} to
lattices $\La$ in elementary locally compact abelian groups because
then the appropriate generalization of $\ZZdba$ is the lattice
$\La^0$ defined as
$\{\la^0\in\RRdd:\pi(\la)\pi(\la^0)=\pi(\la^0)\pi(\la)~~\text{for
all}~~\la\in\La\}$. Feichtinger and Kozek called $\La^0$ the {\it
adjoint lattice} of $\La$. In their discussion of the adjoint
lattice Feichtinger and Kozek pointed out the relevance of the
symplectic Fourier transform in this context. Recently we recognized
that Rieffel had already used the FIGA in his construction of
projective modules over non-commutative tori, \cite{Rief88}. In his
proof of the FIGA he already applied the symplectic Fourier
transform and he defined the adjoint lattice (in his terminology it
is the dual lattice) in full generality $10$ years before
Feichtinger and Kozek. In \cite{FL} we present a comprehensive
treatment of the FIGA and we considerable extend their range of
applications to certain pairs of modulation spaces.
\par
One of the main results in \cite{Jan95} is the following
interpretation of the FIGA for a Gabor system $\cG(g,\al,\be)$.
By definition of $S_{g,\ga,\al,\be}$ the FIGA
expresses that
\begin{equation}
   \langle S_{g,\ga,\al,\be}f,h\rangle=\frac{1}{(\al\be)^d}\langle S_{f,g,\frac{1}{\be},\frac{1}{\al}}\ga,h\rangle
\end{equation}
holds for $f,g,h,\ga\in\SsRd$. More explicitly, the Gabor frame
operator of $\cG(g,\al,\be)$ has the following representation
\begin{equation*}
  S_{g,\ga,\al,\be}f=\frac{1}{(\al\be)^d}\sum_{\klinZd}\langle\ga,\piba g\rangle\piba f.
\end{equation*}
The preceding statement is the so-called {\it Janssen
representation} of a Gabor frame operator. If we start with a Gabor
atom $g\in\LtRd$ then the dual atom $\gamma$ will have the same
quality, i.e. $\gamma\in\LtRd$. Therefore, we have to impose some
extra  conditions to guarantee the validity of the Janssen
representation. These conditions were already introduced by
Tolimieri and Orr in \cite{TO95}. If for a pair $(g,\gamma)$ in
$\LtRd$
\begin{equation*}
  \sum_{\klinZd}\big|\big\langle \gamma,\piba g\big\rangle\big|<\infty,
\end{equation*}
then we say that the pair $(g,\ga)$ satisfies condition $(A')$. If
$g=\ga$, then $g$ is said to satisfy condtion $(A)$. The following
theorem is a consequence of the preceding observations.
\begin{theorem}[Janssen]
Suppose that a pair of functions $(g,\ga)$ in $\LtRd$ satisfies
condition $(A')$ for a given lattice $\ZZab$. Then
\begin{equation*}
  S_{g,\ga,\al,\be}=\frac{1}{(\al\be)^d}\sum_{\klinZd}\big\langle \ga,\piba g\big\rangle\piba
\end{equation*}
holds with absolute convergence in the operator norm.
\end{theorem}

In other words the great insight of Janssen was to express the Gabor
frame operator for nice dual pairs $(g,\ga)$ as an absolutely
convergent sum of time-frequency shifts from the adjoint lattice. We
discuss Janssen's results on the structure of absolutely convergent
sums of time-frequency shifts in section \ref{s:Noncomm} because the
notions of Connes and Rieffel provide the proper framework for these
results.
\par
The result of Tolimieri and Orr that each pair of Schwartz functions
$(g,\ga)$ satisfies the condition $(A')$ led to the search of other
classes of function spaces with this nice property. In addition to the
Schwartz class $\SsRd$ the modulation spaces $\MOvsRd$ provide us with a whole class of functions satisfying the
property $(A')$. We now recall the definition and the basic
properties of $\MOvsRd$.
\par
In $1983$ Feichtinger introduced a class of Banach spaces (see
\cite{Fei83,Fei03}), which allow a measurement of the time-frequency
concentration of a function or distribution $f$ on $\R^d$, the so
called {\sl modulation spaces}. For the measurement of the
time-frequency concentration we choose the {\it Short-Time Fourier
Transform} (STFT) because it is up to a phase-factor the
representation coefficient of the Schr\"odinger representation of
the Heisenberg group. Concretely, if $g$ is a window function in
$\LtRd$ then the STFT of $f\in\LtRd$ is given by
\begin{equation*}
   V_gf(x,\omega):=\langle
   f,\pi(x,\omega)g\rangle=\int_{\R^d}f(t)\ol{g(t-x)}e^{-2\pi i\omega\cdot
   t}dt.
 \end{equation*}

For functions $f$ with good time-frequency concentration, e.g.
Schwartz functions, the STFT can be interpreted as a measure for the
amplitude of the frequency band near $\omega$ at time $x$. The
properties of STFT depend crucially on the window function $g$.

\par
Now we define a function or tempered distribution $f$ to be an
element of the {\it modulation space} $\MOvsRd$ if for a fixed $g$
in Schwartz space $\SsdRd$ the norm
\begin{equation*}
  \|f\|_{\MOvs}:=\|V_gf\|_{\LOvs}=
  \int_{\RRtd}|V_gf(x,\omega)|(1+|x|^2+|\om|^2)^{s/2}dxd\omega
\end{equation*}
is finite. Then $\MOvsRd$ is a Banach space whose definition is
independent of the choice of the window $g$. We always measure the
$\MOvs$-norm with a fixed non-zero window $g\in\SsRtd$. If $s=0$
then modulation space $\MORd$ is the {\it Feichtinger algebra} which
Feichtinger introduced in \cite{Fei81}.
\par
In the following we state some of the properties of $\MOvsRd$, that
are of interest in the later discussion.
  \begin{enumerate}\label{ModspacesThm}

    \item The dual space of $\MOvsRd$ is $\MinfsvRd$ where the duality is given by
    \begin{equation*}
      \langle
      f,h\rangle=\iint_{\R^{2d}}V_{g}f(x,\omega)\ol{V_gh(x,\omega)}dxd\omega,
    \end{equation*}
    for $f\in \MOvsRd$ and $h\in\MinfsvRd$.
    \item \label{TF-invar} $\MOvsRd$ is invariant under time-frequency shifts:
\begin{equation*}
    \|\pi(u,\eta)f\|_{\MOvs}\le
    Cv(u,\eta)\|f\|_{\MOvs} \quad \mbox{for}  \quad  (u,\eta)\in\R^{2d}.
\end{equation*}
    \item $\MOvsRd$ is invariant under Fourier transform.
 \end{enumerate}

\par
The modulation spaces $\MOvsRd$ satisfy a weighted condition $(A')$
which follows from the results in \cite{FG89b}.
\begin{theorem}
If $g,\ga\in\MOvsRd$, then
\begin{equation*}
\sum_{\klinZd}\big|\big\langle \ga,\piba
g\big\rangle\big|(1+\alpha^2|k|^2+\beta^2|l|^2)^{s/2}<C_{\al,\be}\|g\|_{\MOvs}\|\ga\|_{\MOvs}.
\end{equation*}
In particular, if $g,\ga\in\MORd$, then condition $(A')$ is
satisfied simultaneously for all $\al,\be>0$.
\end{theorem}
Finally we mention that the modulation spaces $\MOvsRd$ are the
building blocks of the Schwartz class $\SsRd$:
\begin{equation}
   \SsRd=\bigcap_{s\ge 0}\MOvsRd.
\end{equation}
By duality we get a useful description of the tempered distributions
\begin{equation}
   \SsdRd=\bigcup_{s\ge 0}\MinfsvRd.
\end{equation}
All these statements on modulation spaces are well-known and the
interested reader may find the proofs in Chapter 11 and Chapter 12
of \cite{GrBook}.
\par
We close this section with another important theorem of Janssen. In
his search for nice dual pairs $(g,\ga)$ Janssen took a Schwartz
function $g$ as Gabor atom for the Gabor system $\cG(g,\al,\be)$ and
then he was able to show that the canonical dual window
$\ga_0=S^{-1}_{g,\ga,\al,\be}g$ is also a Schwartz function.

\begin{theorem}[Janssen]
Let $\cG(g,\al,\be)$ be a Gabor frame for $\LtRd$ with $g\in\SsRd$.
Then $\ga_0=S^{-1}_{g,\ga,\al,\be}g\in\SsRd$.
\end{theorem}

In other words if we start with a Gabor frame generated by a
Schwartz function then the canonical dual window has the same
quality. This result has substantially stimulated the research in
Gabor analysis in the last decade. As an outcome of this research
Gr\"ochenig and his collaborators were led to the study of "nice"
classes of involutive Banach algebras \cite{FG05,GL,GL04}.

\section{Basics of non-commutative geometry}\label{s:Noncomm}

\subsection{Motivation}

Around 1980 Connes developed his new idea of non-commutative
differential geometry which generalized the traditional differential
and integral calculus, see \cite{Con94}. Recall that unital and
non-unital commutative $\CS$-algebras correspond to compact and
locally compact spaces, respectively.
\par
One of the central construction in harmonic analysis is the
Pontrjagin dual of a locally compact abelian group. In the case of a
discrete group $\Gamma$ the Pontrjagin dual $\widehat{\Gamma}$ is a
compact abelian group and the duality is given by the Fourier
transform. This shows that the algebra of functions on
$\widehat{\Gamma}$ can be identified with the reduced $C^*$-algebra
of the group $\Gamma$,
\begin{equation*}\label{PontDual}
  C^*_r(\Gamma)\cong C(\widehat{\Gamma}).
\end{equation*}
If $\Gamma$ is non-abelian Pontrjagin duality {\it no longer} works
in the traditional sense, but the left hand side still makes perfect
sense and "behaves" like the algebra of functions on the dual group.
In other words, the Pontrjagin dual $\widehat{\Gamma}$ exists as a
non-commutative space whose algebra of functions is the {\it twisted
group $\CS$-algebra} $C^*_r(\Gamma)$. Therefore non-commutative
geometry provides a natural generalization of Pontrjagin duality, in
the sense that duals of discrete groups are non-commutative spaces.
\par
In this sense unital group $\CS$-algebras provide us with a class of
compact non-commutative spaces, e.g. the group $\CS$-algebra of a
discrete non-abelian group. Consequently one thinks of a non-unital
group $\CS$-algebra as a locally compact non-commutative space. The
group $\CS$-algebra is the most prominent example of this class of
objects. In non-commutative geometry spaces are thought as the same
object if they are Morita equivalent which is a weakening of the
notion of isomorphism for $\CS$-algebras, \cite{Rief72,Rief74a,Rief74b,Rief76}.
\par
Furthermore Connes drew the attention of researchers in operator
algebras to {\it smooth} subalgebras of $C^*$-algebras. The
adjective "smooth" resembles the well-known fact that for a locally
compact Hausdorff space $M$ the space of smooth functions
$C^{\infty}(M)$ is a dense Fr$\acute{e}$chet subalgebra of the
$C^*$-algebra $C(M)$ of continuous functions on $M$. Like
$C^{\infty}(M)$ smooth subalgebras of a $\CS$-algebra $\A$ are in
general Fr$\acute{e}$chet subalgebras of $\A$. In the last two
decades "smooth subalgebras" of $C^*$-algebras have been studied by
many researchers, e.g. Badea, Bost, Cuntz, Ji and Schweitzer, see
\cite{Sch92,Sch93,Sch94}.
\par
In \cite{GL,GL04} Gr\"ochenig and Leinert have concentrated their
efforts on non-commutative Banach algebras of a $\CS$-algebra $\A$
which may considered as the non-commutative analog of the $k$-times
differentiable functions $C^r(M)$ on $M$. Later we expose their
results for rotation algebras in detail. After this short digress on
the underlying principles of non-commutative geometry we introduce
the main operator algebras of our investigations: (1) rotation
algebras and (2) non-commutative tori.

\subsection{Non-commutative tori}
Since the 1960's rotation algebras have been studied by Effros and
his collaborators but the work \cite{Con80,Rief81} of Connes and
Rieffel on the structure of rotation algebras brought these algebras
into the focus of researchers outside of operator algebras.
\par
It is well-known that the $\CS$-algebra of continuous functions on
the torus is the universal $\CS$-algebra generated by two commuting
unitaries, which can be considered as the coordinate functions. This
suggests considering the universal $\CS$-algebra $\TT_\theta^2$
generated by unitaries $U,V$ satisfying the commutation relation
\begin{equation*}
  UV=e^{2\pi i\theta}VU,
\end{equation*}
where the $\theta$ is a real number. $\TT_\theta^2$ is called the
{\it rotation algebra}. For an insightful exposition of rotation
algebras we refer the reader to the monograph \cite{Dav96}. We
restrict ourselves to the connection between Gabor analysis and
rotation algebras.
\par
 Let $\theta=\alpha\beta$ for positive reals $\alpha$ and $\beta$ then the
$C^*$-algebra generated by time-frequency shifts $\{\pi(\alpha
k,\beta l) :k,l\in\ZZ^{d}\}$ is a representation of the rotation
algebra $\TT_{\theta}^2$ on $L^2(\RR^d)$. Therefore, an element of
$\TT_\theta^2$ is given by
\begin{equation*}
  \sum_{\klinZd}a_{kl}\piab
\end{equation*}
for a bounded complex-valued sequence ${\bf a}=(a_{kl})_{\klinZd}$.
\begin{lemma}
The representation of the rotation algebra by time-frequency shifts
$\{\piab:\klinZd\}$ is {\bf faithful} on $L^2(\RR^d)$.
\end{lemma}
As a consequence it suffices to establish statements about the
rotation algebra $\TT_\theta^2$ for a dense subspace of
$L^2(\RR^d)$. In \cite{Rief88} the interested reader finds a proof
based on operator algebra methods and in \cite{GL04} the reader finds
another proof relying on time-frequency methods and Wiener amalgam
spaces.
\par
The faithfulness of the representation of $\TT_\theta^2$ by
time-frequency shifts $\piab$ from the lattice $\ZZdab$ is essential
for the results of Janssen, Gr\"ochenig-Leinert and Rieffel,
\cite{GL04,Jan95,Rief88}.
\par
In \cite{Con80} a {\it non-commutative torus} $\A^{\infty}$ was
defined as the "smooth" elements of $\TT_\theta^2$, i.e.
\begin{equation*}
  \A^{\infty}=\{A\in\mathcal{B}\big(L^2(\RRd)\big):A=\sum_{k,l}a_{kl}U^kV^l,~~{\bf
  a}=(a_{kl})_{\klinZd}\in\Ss(\ZZtd)\}.
\end{equation*}
$\Ss(\ZZtd)$ is the space of sequences with rapid decay on $\ZZtd$.
$\A^{\infty}$ was the first example of a {\it smooth} structure on a
non-commutative $\CS$-algebra.
\par
In 1995 Janssen described the non-commutative torus $\A^{\infty}$ in
the representation of $\TT_\theta^2$ by time-frequency shifts $\piab$.
More precisely, in \cite{Jan95} $\A^{\infty}$ is described as the
following family of Banach algebras:
\begin{equation*}
  \A^{1}_s(\al,\be)=\{A\in\mathcal{B}\big(L^2(\RRd)\big):
  A=\sum_{\klinZd}a_{kl}\piab,~~{\bf
  a}=(a_{kl})_{\klinZd}\in\ell^1_s(\ZZtd)\},
\end{equation*}
where $\ell^1_s(\ZZtd)$ is the space of all sequences on $\ZZtd$
such that
\begin{equation*}
\|{\bf
a}\|_{1,s}:=\sum_{\klinZd}|a_{kl}|(1+\al^2|k|^2+\beta^2|l^2|)^{s/2}<\infty.
\end{equation*}
Consequently, the non-commutative tori
$\A^{\infty}(\al,\be)$ is just $\bigcap_{s\ge 0}\A^{1}_s(\al,\be)$.
\par
The Banach algebras $\ell^{1}_s(\ZZtd)$ inherit from the commutation
relations
\begin{equation*}
  \piab\pi(\al k',\be l')=e^{2\pi i(k'\cdot l-k\cdot l')}\pi(\al k',\be l')\piab
\end{equation*}
a {\it twisted product}
\begin{equation*}
   \big({\bf a}\natural{\bf b}\big)(m,n)=\sum_{\klinZd}a_{kl}b_{m-k,n-l}e^{2\pi i\theta(m-k)\cdot l}
\end{equation*}
for ${\bf a}, {\bf b}\in\ell^{1}_s(\ZZtd)$ and a {\it twisted
involution}
\begin{equation*}
   a^{\ast}(k,l)=\ol{a_{kl}}e^{-2\pi i\theta k\cdot l}
\end{equation*}
for ${\bf a}=(a_{kl})_{\klinZd}\ell^{1}_s(\ZZtd)$.
\begin{lemma}
$\big(\A^{1}_s(\al,\be),\natural,\ast\big)$ is an involutive Banach
algebra.
\end{lemma}

\par
We consider the involutive Banach algebra $\A^1_s(\ZZt)$ as the
non-commutative analog of $C^s$-functions for the rotation algebra
$\TT_\theta^2$. Since the non-commutative analog of the
differentiation operator is a {\it derivation} $\delta$ of a
$\CS$-algebra $\A$, i.e.
\begin{equation*}
  \delta(AB)=\delta(A)B+A\delta(B),~~A,B\in\A.
\end{equation*}
On the non-commutative torus $\A^{\infty}$ we have a pair of
commuting derivations
\begin{equation*}
  \delta_1\big(\sum a_{kl}U^kV^l\big)=2\pi i \sum ka_{kl}U^kV^l
\end{equation*}
and
\begin{equation*}
 \delta_2\big(\sum a_{kl}U^kV^l\big)=2\pi i \sum la_{kl}U^kV^l
\end{equation*}
for $A=\sum a_{kl}U^kV^l$. Therefore the non-commutative analog of
the Laplacian $\Delta=\delta_1^2+\delta_2^2$ acts on $\A^{\infty}$
by
\begin{equation*}
 \big(\delta_1^2+\delta_2^2\big)\big(\sum a_{kl}U^kV^l\big)=-4\pi^2\sum (k^2+l^2)a_{kl}U^kV^l.
\end{equation*}
By this reasoning we are tempted to consider those
$A\in\TT_\theta^2$ such that $\|\Delta^{s/2} A\|_{\text{op}}<\infty$
as non-commutative analog of $C^s$ for $s\in[0,\infty)$. Now, we
suitably normalize our "Laplacian" $\tilde\Delta$ and define the
non-commutative analog of the potential operator by $I-\tilde\Delta$
which acts on $A\in\A^{\infty}$ in the following way,
\begin{equation*}
 \big(I-\til\Delta\big)\big(\sum a_{kl}U^kV^l\big)=\sum (1+k^2+l^2)a_{kl}U^kV^l.
\end{equation*}
The preceding discussion suggests to think of the space of elements
$A\in\TT_\theta^2$ such that
$\|(I-\til\Delta)^{s/2})A\|_{\text{op}}<\infty$ for $s\in[0,\infty)$
as the non-commutative analog of the Sobolev spaces. In the case of
the representation of  $\TT_\theta^2$ by $\{\piab:\klinZd\}$ the
space of $\|(I-\til\Delta)^{s/2})A\|_{\text{op}}<\infty$ turns out
to be $\A^1_s(\al,\be)$.
\par
Recall the fact that $f,g\in\MOvsRd$ implies that
$\big\|\big(\langle f,\piab g\rangle\big)_{k,l}\big\|_{1,s}<\infty$,
i.e.
\begin{equation*}
  A=\sum_{\klinZd}\langle f,\piab g\rangle\piab\in\\A^1_s(\al,\be).
\end{equation*}
The preceding observation indicates the relevance of modulation
spaces $\MOvsRd$ for Gabor analysis and non-commutative tori. We
close this section with a characterization of $\MOvsRd$ by means of
Gabor frames.
\begin{theorem}
  Assume that $g,\ga\in\MOvsRd$ and that $S_{g,\gamma,\al,\be}=I$ on $\LtRd$. Then
  \begin{equation*}
    f=\sum_{\klinZd}\langle f,\piab g\rangle\piab\gamma=\sum_{\klinZd}\langle f,\piab \gamma\rangle\piab g
  \end{equation*}
  with unconditional convergence in $\MOvsRd$.
  \par Furthermore, there are constants $A,A',B,B'>0$ such that for all $f\in\MOvsRd$
  \begin{equation*}
A\|f\|_{\MOvs}\le\sum_{\klinZd}|\langle f,\piab g\rangle|\le
B\|f\|_{\MOvs}
  \end{equation*}
  and
\begin{equation*}
A'\|f\|_{\MOvs}\le\sum_{\klinZd}|\langle f,\piab \gamma\rangle|\le
B'\|f\|_{\MOvs}.
  \end{equation*}
\end{theorem}
The preceding result is just a special case of a general principle
but it provides another reason for the search of nice classes of windows, see
\cite{GrBook}.

\section{Spectral invariance of Banach and Frechet algebras}

In this section we proceed to state the non-commutative analogs of
some traditional notions. In Connes' calculus a complex variable
corresponds to an operator $T$ on an infinite-dimensional Hilbert
space $\cH$ and a real variable is associated with a self-adjoint
operator on $\cH$. Consequently, the {\it spectrum} $\sigma(T)$ of
$T$ is the non-commutative analog of the range of a complex
variable. The holomorphic functional calculus for operators in a
Hilbert space gives meaning to  $f(T)$ for any holomorphic function
$f$ defined on $\sigma(T)$ and only holomorphic functions act in
that generality. This reflects the need for holomorphy in the theory
of complex variables. Indeed, when the operator $T$ is self-adjoint,
$f(T)$ makes sense for any Borel function $f$ on the line. In the
following we present an extensive discussion of the holomorphic
functional calculus. We begin with some well-known results of
Wiener, Gelfand and Shilov on absolutely convergent Fourier series.
\par
If we consider the Gabor system $\cG(g,1,1)$ then the Banach
algebras $\A^1_s(\al,\be)$ and the Fr$\acute{e}$chet algebra
$\A^{\infty}(\al,\be)$ are commutative and involutive subalgebras of
$\TT^2$, i.e. the spaces of functions with absolutely convergent and rapidly decaying Fourier series. 
Now, an application of the Fourier transform allows the
inversion of the Gabor frame operator, \cite{FG97}.
\par
More concretely, Let ${\mathcal A}(\TT^2)$ be the Banach algebra of all absolutely
convergent Fourier series, i.e.
\begin{equation*}
  {\mathcal A}(\TT^2)=\{f\in C(\TT^2):f(x,t)=\sum_{k,l\in\ZZd}a_{kl}e^{2\pi i(kx+lt)},~~\sum_{k,l\in\ZZd}|a_{kl}|<\infty\}
\end{equation*}
with norm $\|f\|_{{\mathcal A}(\TT^2)}=\sum_{k,l\in\ZZd}|a_{kl}|$.
In the early $30$'s of the last century N. Wiener made an important
observation about non-zero functions in ${\mathcal A}(\TT^2)$.
Namely, if $f\in{\mathcal A}(\TT^2)$ is non-zero on $\TT^2$, then
$1/f$ is in ${\mathcal A}(\TT^2)$, too. More precisely, there exists
a sequence $b=(b_{kl})_{k,l\in\ZZd}$ such that
$1/f=\sum_{k\in\ZZd}b_{kl}e^{2\pi i(kx+lt)}$. This fact goes by the
name {\it Wiener's lemma}. This result of Wiener has many important
applications and therefore ${\mathcal A}(\TT^2)$ is often referred
to as {\it Wiener's algebra}.
\par
The first great success of Gelfand's theory of Banach algebras was a
short elusive proof of Wiener's lemma. Later Naimark pointed out
that Wiener's lemma is actually a statement about the pair of
involutive Banach algebras ${\mathcal A}(\TT^2)\subset C(\TT^2)$
such that if a function $f\in{\mathcal A}(\TT^2)$ is invertible in
$C(\TT^2)$ then $1/f$ is an element of ${\mathcal A}(\TT^2)$. Naimark
inspired by Wiener's result introduced the following notion.
\begin{definition}
Let ${\mathcal A}$ be a subalgebra of the Banach algebra ${\mathcal
B}$ with common unit $I$. Then ${\mathcal A}$ is called {\it
spectral invariant} in ${\mathcal B}$ if $A\in{\mathcal A}$ and
$A^{-1}\in{\mathcal B}$ implies $A^{-1}\in{\mathcal A}$. In this
case $({\mathcal A},{\mathcal B})$ is called a {\it Wiener pair}.
\end{definition}
Gr\"ochenig and his collaborators call $\A$ {\it inverse-closed} in $\B$.
\par
At the moment we draw some general conclusions from the definition
of spectral invariance. We will focus on the spectral radius, the
spectrum and holomorphic functional calculus for a general Wiener
pair $({\mathcal A},{\mathcal B})$ of Banach algebras. Since $A$ is
an element of ${\mathcal A}$ and ${\mathcal B}$ we can talk about
the spectrum of $A$ with respect to the algebra ${\mathcal A}$ and
${\mathcal B}$, respectively. Recall that the spectrum of $A$ in
${\mathcal A}$ is defined as
\begin{equation}
  {\sigma}_{{\mathcal A}}(A)=\{z\in\CC:A-zI~~\text{is not invertible in}~~{\mathcal A}\}.
\end{equation}
 An elementary argument yields the following lemma.
\begin{lemma}
For ${\mathcal A}\subset{\mathcal B}$ with common unit $I$ then
\begin{equation*}
  \sigma_{{\mathcal B}}(A)\subset\sigma_{{\mathcal A}}(A).
\end{equation*}
\end{lemma}
In combination with the preceding observation we get the following
innocent looking lemma which gives a justification to call a Wiener pair of Banach algebras $(\A,\B)$ spectral invariant.
\begin{lemma}
  Let ${\mathcal A}\subset{\mathcal B}$ be Banach algebras with common unit $I$. Then the following
statements are equivalent:
\begin{enumerate}
  \item $({\mathcal A},{\mathcal B})$ is a Wiener pair.
  \item $\sigma_{{\mathcal A}}(A)=\sigma_{{\mathcal B}}(A).$
\end{enumerate}
\end{lemma}
\par
In many situations we have to invert or to take a square root of an
element $A$ in a Banach algebra ${\mathcal A}$ with unit element
$I$. In the early years F. Riesz, one of the pioneers of functional
analysis, had the idea to build functions $f(T)$ of an compact
operator $T\in{\mathcal B}(\LtR)$ in analogy to Cauchy's formula in
complex analysis. Later this kind of reasoning had been continued by
Wiener and substantially generalized by Dunford. Therefore the
calculus goes by the names: {\it Riesz functional calculus}, {\it
Dunford calculus} or {\it holomorphic calculus}.
\par
More precisely, let ${\mathcal A}$ be a unital Banach algebra and
$A\in{\mathcal A}$. Then
\begin{equation}
  \text{Hol}(A)=\{f:f~\text{is holomorphic on an open neighborhood}~G~\text{of}~\sigma_{{\mathcal A}}(A)\}
\end{equation}
is the reservoir of functions which allows to form new elements
$\tilde f(A)$ in the Banach algebra ${\mathcal A}$. Therefore, we
choose a neighborhood $G$ of $\sigma_{{\mathcal A}}(A)$ and a
contour $\Gamma$ of $\sigma_{{\mathcal A}}(A)$ in $G$. Then
$f\in\text{Hol}(A)$ with domain $G$ allows to define
\begin{equation}
  \tilde f(A):=\frac{1}{2\pi i}\int_{\Gamma}f(z)(zI-A)^{-1}dz
\end{equation}
as a Banach algebra valued integral. The basic results on Cauchy's
formula in complex analysis imply that the definition of $\tilde
f(A)$ is independent of $\Gamma$ and that the integral exists as a
Riemann integral and $\tilde f(A)$ is an element of ${\mathcal A}$.
Recall, that $R(z,A)=(zI-A)^{-1}$ is called the {\it resolvent
function} of $A$ which is defined on the {\it resolvent set}
$\rho_{\mathcal A}(A)=\CC\backslash\sigma_{\mathcal A}(A)$. The
resolvent function $R(.,A)$ is analytic on $\rho_{\mathcal A}(A)$
and therefore $z\mapsto f(z)(zI-A)^{-1}$ is analytic from
$G\cap\rho_{\mathcal A}(A)$ into ${\mathcal A}$.
\par
The main theorem in this context says that for a fixed
$A\in{\mathcal A}$. The mapping $A\mapsto\tilde f(A)$ is an algebra
homomorphism and that this mapping is continuous from
$\text{Hol}(A)$ under uniform convergence on compact sets to
${\mathcal A}$ with the norm topology.
\par
The following result is what one needs in the discussion of "nice"
window classes of Gabor frames.
\begin{theorem}
Let ${\mathcal A}\subset{\mathcal B}$ be Banach algebras with common
unit $I$. If $({\mathcal A},{\mathcal B})$ is a Wiener pair, then
the Riesz functional calculus for ${\mathcal A}$ coincides with the
one for ${\mathcal B}$.
\end{theorem}
\begin{proof}
  Since $\sigma_\A(A)=\sigma_\B(A)$, the resolvent function $R(z,.)$ is
  defined in $\A$ and $\B$. Consequently, we get that $\text{Hol}_\A(A)$ and
  $\text{Hol}_\B(A)$ coincide.
\end{proof}

The preceding discussion tells us that a Wiener pair $({\mathcal
A},{\mathcal B})$ has many nice properties. But how can we decide if
two Banach algebras ${\mathcal A}\subset{\mathcal B}$ with common
unit $I$ form a Wiener pair? That's actually the hard part in this
topic and the only known tool is the so-called {\it lemma of
Hulanicki} which relates the spectral invariance of $({\mathcal
A},{\mathcal B})$ with the symmetry of ${\mathcal B}$.
\par
The notion of symmetry of an involutive Banach algebra ${\mathcal
A}$ with unit $I$ is the generalization of the fact that a positive
bounded operator $T$ on a Hilbert space ${\mathcal H}$ has spectrum
$\sigma(T)\subset [0,\infty)$.

\begin{definition}
An involutive Banach algebra ${\mathcal A}$ with unit $I$ is called
{\bf symmetric}, if $\sigma_{\mathcal A}(AA^*)\subset[0,\infty)$ for all $A\in{\mathcal A}$.
\end{definition}
An element $A\in{\mathcal A}$ is called {\bf positive} if $A=CC^*$
for some $C\in{\mathcal A}$. In this terminology we can say that
positive elements $A$ in an involutive Banach algebra ${\mathcal A}$
with unit have "positive" spectrum. Therefore symmetry of an
involutive Banach algebra measures how close ${\mathcal A}$ is to be
a $C^*$-algebra.
\begin{remark}
An involutive Banach algebra ${\mathcal A}$ with unit is symmetric
if and only if $\sigma_{\mathcal A}(A)\subset\RR$ for all $A=A^*$ in
${\mathcal A}$, if $A=A^*$ then $A$ is {\bf hermitian}. This fact is
the {\bf Ford-Shiraly lemma} and has been an open question for many
years in the early days of normed algebras.
\end{remark}
The following theorem is the lemma of Hulanicki from the early
$70$'s.
\begin{theorem}[Lemma of Hulanicki]
Assume that ${\mathcal B}$ is a symmetric Banach algebra and
${\mathcal A}$ a subalgebra of ${\mathcal B}$ with common unit $I$.
Then ${\mathcal A}$ is spectral invariant in ${\mathcal B}$ if and
only if the spectral radii for all $A=A^*$ with respect to
${\mathcal A}$ and ${\mathcal B}$ are equal, i.e.
$\operatorname{spr}_{\mathcal A}(A)=\operatorname{spr}_{\mathcal
B}(A)$
\end{theorem}
A recent unpublished result by Leinert relates the notions of
spectral invariance and symmetry of an involutive Banach algebra $\A$.
\begin{theorem}[Leinert]\label{Leinert}
Let $C^*({\mathcal A})$ be the enveloping $C^*$-algebra of an
involutive Banach algebra ${\mathcal A}$ with unit $I$ then
${\mathcal A}$ is symmetric if and only if ${\mathcal A}$ is
spectral invariant in $C^*({\mathcal A})$.
\end{theorem}
\par

Now we have all notions and tools in hand to continue our
investigations of "nice" window classes of a Gabor frame operator.
Before we proceed we have to make some remarks. In \cite{Jan95} the
result that the canonical dual atom has the same quality as the
Gabor atom was only proved under the additional assumption that
$\theta=\al\be$ is a rational number. At the end of the paper
Janssen formulated the conjecture that the result is also true if
$\theta$ is irrational. Therefore this conjecture was called the
{\it irrational case conjecture}. The resolution of Janssen's
conjecture was the great stimulus for the work of Gr\"ochenig and
Leinert on non-commutative analogs of Wiener's lemma for the
irrational rotation algebra, \cite{GL04}.

\begin{theorem}[Gr\"ochenig-Leinert]\label{GrLe}
 ${\mathcal A}^1_s(\al,\be)$ is spectral invariant in $\TT_\theta^2$.
\end{theorem}
\begin{corollary}
${\mathcal A}^1_s(\al,\be)$ is stable under holomorphic calculus.
\end{corollary}
If one recalls that the rotation algebra $\TT_\theta^2$ is the
enveloping $\CS$-algebra of $\ZZtd$ for the cocyle
$\chi((k,l),(m,n))=e^{2\pi im\cdot l}$. In \cite{GL04} the main
result is the symmetry of the twisted convolution algebra
$\ell^1_s(\ZZtd,\natural,^\ast)$.
\par
Consequently, the theorem of Gr\"ochenig and Leinert is a special case of Theorem \ref{Leinert}.
At the beginning of the section we discussed Wiener's lemma and
Wiener's algebra $\A(\TT^2)$. Therefore the Banach algebra
$\A^1_s(\al,\be)$ is called the {\it non-commutative Wiener
algebra}. As a further corollary of Theorem \ref{GrLe} we state the
following result.
\begin{theorem}
Let $\ZZdab$ a lattice in $\RRtd$ and let $g\in\MOvsRd$. If
$\cG(g,\al,\be)$ is a Gabor frame for $\LtRd$, then
$S_{g,\al,\be}^{\nu}g\in\MOvsRd$ for $\nu\in\RR$. Especially, we
have that for any $f\in\LtRd$
\begin{equation*}
       f=\sum_{\klinZd}\langle f,\piab g\rangle\piab S_{g,\al,\be}^{-1}g
         =\sum_{\klinZd}\langle f,\piab S_{g,\al,\be}^{-1/2}g\rangle\piab S_{g,\al,\be}^{-1/2}g.
\end{equation*}
\end{theorem}
Another consequence of Theorem \ref{GrLe} is the spectral invariance
of $\A^{\infty}(\al,\be)$. Since our weights $v_s(x,\om)$ are
submultiplicative for $s\in\RR$, i.e.
$v_s(x+y,\om+\eta)\le v_s(x,\om)v_s(y,\eta)$. Therefore
$\A^{\infty}(\al,\be)$ is a Fr$\acute{e}$chet algebra with a scale
of Banach algebras $\A^1_s(,\al,\be)$ determined by a
submultiplicative seminorm. Now we invoke an old result from the
early days of operator algebras due to Michael, which says that the
spectral invariance of each level $\A^{1}_s(\al,\be)$ implies the
spectral invariance of $\A^{\infty}(\al,\be)$, \cite{Mic52}.
\begin{corollary}[Connes-Janssen]
$\A^{\infty}(\al,\be)$ is spectral invariant in the rotation algebra
$\TT_\theta^2$.
\end{corollary}
We close the section with a short review of Connes' argument for the spectral invariance of $\A^{\infty}$ and $\TT_\theta^2$. 
\par
We call a subalgebra $\A$ of a unital Banach algebra $\B$ {\bf
closed under holomorphic functional calculus} of $ A$ if it
satisfies the following conditions:
\begin{enumerate}
  \item $\A$ is complete under some locally convex topology finer than the topology of $\B$;
  \item If $A\in\A$ and $f(A)$ is defined by the Riesz-Dunford integral, then $f(A)\in\A$.
\end{enumerate}

In \cite{Con80} Connes coined the notion of a {\it
pre-$\CS$-algebra} for a subalgebra $\A$ of a $\CS$-algebra $\B$
that is closed under holomorphic functional calculus.
\par
The most well-known example of a pre-$\CS$-algebra is the Schwartz
class ${\mathcal S}(\RRd)$ whose $\CS$-completion is the
$\CS$-algebra of continuous functions $C_0(\RRd)$ on $\RRd$. The
characterization of the Schwartz class $\Ss(\RRd)=\cap_{s\ge
0}M^1_{v_s}(\RRd)$. In non-commutative geometry all pre-$\CS$-algebras are
Fr$\acute{e}$chet algebras like $C^\infty(\RRd)$ and $\Ss(\RRd)$. Furthermore they 
arise as smooth vectors of a strongly continuous action.
\begin{lemma}
  Let $\A$ be a $\CS$-algebra, $G$ a Lie group and $\rho$ a strongly continuous action from $G$ into the automorphisms $\text{Aut}(\A)$ of $\A$. The dense subalgebra $\A^\infty$ of smooth elements under this action is a Fr$\acute{e}$chet pre-$\CS$-algebra.
\end{lemma}
\begin{proof}
  Now, for each $A\in\A$ the map $t\mapsto\rho_t(A)$ is continuous and $\A^\infty$ consists of those $A\in\A$ such that this map is smooth. Since for any $f\in\text{Hol}(A)$ we have $f(\rho_t(A))=\rho_t(f(A))$, the algebra $\A^\infty$ is closed under holomorphic functional calculus in $\A$.
\end{proof}
Now, Connes considered the following strongly continuous action. Let
$\TT^{2d}$ be the dual group of $\ZZtd$. Then $\TT^{2d}$ has a
natural dual action $\rho$ on $\TT^2_\theta$ given by
\begin{equation*}
  (\rho_t)({\bf a})(k)=e^{2\pi i t\cdot k}a(k),~~\text{for}~~{\bf a}\in\ell^1_s(\ZZtd),k\in\ZZtd,t\in\TT^{2d}.
\end{equation*}
By the first lemma of Section 13 of \cite{Con82}, the space of
smooth vectors for this action will be exactly $\Ss(\ZZtd)$. Now
Connes used that $\Ss(\ZZtd)$ is closed under holomorphic functional
calculus, see the appendix in \cite{Con80}. Therefore
$\A^{\infty}(\al,\be)$ is spectral-invariant in $\TT^2_\theta$.

\section{Conclusion}

The present investigation is aimed to link Janssen's outstanding
contribution to Gabor analysis with non-commutative geometry,
especially Connes's discussion of non-commutative tori. Our
intention was to stress that a problem in applied mathematics gives
rise to a deep result on non-commutative tori. More concretely,
Gr\"ochenig and Leinert's main insight is the equivalence of nice classes of 
Gabor atoms with the study of spectral invariant subalgebras of rotation algebras. 
But they were not aware of Connes's contribution to this problem and its connection 
to non-commutative geometry. Furthermore Gr\"ochenig and Leinert restricted themselves to spectral invariant Banach
algebras contrary to Connes's search for spectral invariant
Fr$\acute{e}$chet subalgebras of rotation algebras. 
Finally we want to mention that Janssen's results have allowed to design new algorithms in
Gabor analysis which lead to a significiant improvement of
transfer rates for cellular phones or OFDM networks, see the
contributions by Kozek and Strohmer in \cite{FS98}. In other words
{\it non-commutative tori have real-world applications}.

\section{Acknowledgment}
The author wants to thank K. Gr\"ochenig for several discussions on
spectral invariant Banach algebras and for his inspiring course
"Banach algebra methods in applied mathematics" held at the
University of Vienna during spring term 2005. Furthermore I express
my gratitude for the kind invitation of D. Han and D. Larson to
present some results of my Ph.D thesis at GPOTS 2005. Finally I want to thank 
H.G. Feichtinger for many suggestions and comments which improved the presentation considerably.

\flushleft{
\bibliographystyle{amsalpha}

\end{document}